\documentclass{amsart}
\usepackage{amsfonts,amssymb,epsfig}
\usepackage[latin1]{inputenc}
\newtheorem{thm}{Theorem}[section]

\newtheorem{lemma}[thm]{Lemma}
\newtheorem{prop}[thm]{Proposition}
\theoremstyle{definition}

\theoremstyle{remark}

\numberwithin{equation}{section}

\renewcommand{\epsilon}{\varepsilon}

\newcommand{\fg}{\mathfrak{g}}

\newcommand{\cY}{\mathcal Y}

\newcommand{\cO}{\mathcal O}
\newcommand{\cR}{\mathcal R}
\newcommand{\bbR}{\mathbb R}

\newcommand{\bbC}{\mathbb C}
\newcommand{\bbK}{\mathbb K}
\newcommand{\bbN}{\mathbb N}

\begin{document}

\title[Levi decomposition of Poisson structures and Lie algebroids]
{Levi decomposition of analytic Poisson structures and Lie algebroids}

\author{Nguyen Tien Zung}
\address{Laboratoire Emile Picard, UMR 5580 CNRS, UFR MIG, Université Toulouse III}
\email{tienzung@picard.ups-tlse.fr} 

% 53D17 : Poisson manifolds
% 32S65 Singularities of holomorphic vector fields and foliations
\date{second version, Dec/2002, to appear in Topology}%

\begin{abstract}
We prove the existence of a local analytic Levi decomposition for analytic Poisson
structures and Lie algebroids.
\end{abstract}
\maketitle

\section{Introduction}

In the study of local normal forms of Poisson structures, initiated by Arnold
\cite{Arnold-Small1963} Weinstein \cite{Weinstein-local}, one is led naturally to
the following problem of Levi decomposition: let $\Pi$ be a Poisson structure in a
neighborhood of $0$ in $\bbK^n$, where $\bbK = \bbR$ or $\bbC$, such that $\Pi(0) =
0$. We will use the letter $\Pi$ to denote the Poisson tensor, and $\{.,.\}$ or
$\{.,.\}_\Pi$ to denote the corresponding Poisson bracket. In this paper we will
assume that $\Pi$ is analytic. Denote by $\Pi_1$ the linear part of $\Pi$ at $0$.
$\Pi_1$ is a linear Poisson tensor, and the space $\frak L$ of linear functions on
$\bbK^n$ is an $n$-dimensional Lie algebra under the Poisson bracket of $\Pi_1$.
Denote by $\frak r$ the radical of $\frak L$. The classical Levi-Malcev theorem says
that the exact sequence $0 \to \frak r \to \frak L \to {\frak L}/{\frak r} \to 0$
admits a splitting : there is an injective homomorphism from  ${\frak L}/{\frak r}$
to $\frak L$ (unique up to a conjugation in $\frak L$) whose composition with the
projection map is identity. Denote by $\fg$ the image of such an inclusion. Then
$\fg$ is called a Levi factor of $\frak L$, and $\frak L$ can be written as a
semi-direct product of a semi-simple Lie algebra $\fg$ by a solvable Lie algebra
$\frak r$ (this semi-direct product is called a Levi decomposition of $\frak L$).
Remark that the space $\cO$ of local analytic functions in $(\bbK^n,0)$ is an
infinite-dimensional Lie algebra under the Poisson bracket $\{.,.\}_\Pi$, and  the
space $\cR$ of local analytic functions in $(\bbK^n,0)$ whose linear part lies in
$\frak r$ is an infinite-dimensional ``radical'' of $\mathcal O$, with $\cO / \cR$
isomorphic to $\fg$. The question is, does the exact sequence $0 \to \cR \to \cO \to
\fg \to 0$ also admit a splitting ? In other words, does $\cO$ together with the
Poisson structure $\Pi$ admit a Levi factor ? In this paper, we will give a positive
answer to this question. More precisely, we have :

\begin{thm}
\label{thm:LeviPoisson} Let $\Pi$ be a local analytic Poisson tensor in
$(\bbK^n,0)$, where $\bbK = \bbR$ or $\bbC$. Denote by $\frak L$ the $n$-dimensional
Lie algebra of linear functions in $(\bbK^n,0)$ under the Lie-Poisson bracket of
$\Pi_1$ which is the linear part of $\Pi$, and by $\frak L = \fg \ltimes \frak r$ a
Levi decomposition of $\frak L$. Denote by $(x_1,...,x_m,y_1,...,y_{n-m})$ a linear
basis of $\frak L$, such that $x_1,...,x_m$ span the Levi factor $\fg$ ($\dim \fg =
m$), and $y_1,...,y_{n-m}$ span the radical $\frak r$. Denote by $c_{ij}^k,
b_{ij}^k$ and $a_{ij}^k$ the structural constants of $\fg, \frak r$ and of the
action of $\fg$ on $\frak r$ respectively : $[x_i,x_j] = \sum_k c_{ij}^k x_k$,
$[y_i,y_j] = \sum_k b_{ij}^k y_k$ and $[x_i,y_j] = \sum_k a_{ij}^k y_k$. Then there
exists a local analytic system of coordinates
$(x^{\infty}_1,...,x^{\infty}_m,y^{\infty}_1,..., y^{\infty}_{1-m})$, with
$x^{\infty}_{i} = x_i +$ higher order terms and $y^{\infty}_{i} = y_i +$ higher
order terms, such that in this system of coordinates we have
\begin{equation}
\label{eqn:LeviPoisson1} \Pi = \frac{1}{2} [\sum c_{ij}^k x^{\infty}_{k}
\frac{\partial}{\partial x^{\infty}_{i}} \wedge \frac{\partial }{ \partial
x^{\infty}_{j}} + \sum a_{ij}^k y^{\infty}_{k} \frac{\partial }{
\partial x^{\infty}_{i}} \wedge \frac{\partial }{ \partial
y^{\infty}_{j}} + \sum (b_{ij}^k y^{\infty}_{k} + g_{ij}) \frac{\partial }{
\partial y^{\infty}_{i}} \wedge \frac{\partial }{ \partial y^{\infty}_{j}}]
\end{equation}
where $g_{ij}$ are local analytic functions whose Taylor expansion begins at order
at least 2. In other words, the Poisson bracket $\{.,.\}$ of $\Pi$ in this system of
coordinates is given as follows :
\begin{equation}
\label{eqn:LeviPoisson2}
\begin{array}{l}
\{x^{\infty}_{i},x^{\infty}_{j}\} = \sum c_{ij}^k x^{\infty}_{k} \ , \\
\{x^{\infty}_{i},y^{\infty}_{j}\} = \sum a_{ij}^k y^{\infty}_{k}  \ , \\
\{y^{\infty}_{i},y^{\infty}_{j}\} = \sum b_{ij}^k y^{\infty}_{k} + g_{ij} \ .
\end{array}
\end{equation}

\end{thm}

{\it Remarks.}

1. In the above theorem, the Levi factor of $\cO$  is provided by the functions
$x^{\infty}_{1},...,x^{\infty}_m$. Conversely, if $\cO$ admits a Levi factor with
respect to $\Pi$, then the Hamiltonian vector fields of the functions lying in this
Levi factor gives us a local analytic Hamiltonian action of $\fg$, which is
linearizable by a theorem of Guillemin and Sternberg \cite{GuSt-Linearization1968},
because $\fg$ is semi-simple. By linearizing this action, one will get a local
analytic coordinate system which satisfies the conditions of the above theorem. Thus
the above theorem is really about the existence of an analytic Levi decomposition of
the Poisson structure.

2. If in the above theorem, we don't require the functions
$x^{\infty}_1,...,x^{\infty}_m,y^{\infty}_1,..., y^{\infty}_{1-m}$ to be analytic,
but only formal, then we recover a formal Levi decomposition theorem, obtained
earlier by Wade \cite{Wade}. This formal decomposition is relatively simple and its
proof is similar to the proof of the classical Levi-Malcev theorem. The difficulty
of the above theorem lies in the analytic part.

3. If in the above theorem, $(\frak L, \{.,.\}_{\Pi_1})$ is a semi-simple Lie
algebra, i.e. $\fg = \frak L$, then we recover the following result of Conn
\cite{Conn} : any analytic Poisson structure with a semi-simple linear part is
locally analytically linearizable. In other words, any semi-simple Lie algebra is
analytically nondegenerate in the terminology of Weinstein \cite{Weinstein-local}.
In fact, our proof of Theorem \ref{thm:LeviPoisson} will follow closely the lines of
Conn \cite{Conn}. When $\frak r = \bbK$ ($\bbK = \bbR$ in the real case and $\bbK =
\bbC$ in the complex case), i.e., $\frak L = \fg \oplus \bbK $, we get the following
result, due to Molinier \cite{Molinier} and Conn (unpublished) : if $\fg$ is
semi-simple then $\fg \oplus \bbK$ is analytically nondegenerate.

4. One may call expressions (\ref{eqn:LeviPoisson1}), (\ref{eqn:LeviPoisson2}) a
Levi normal form of the Poisson structure $\Pi$. From the point of view of invariant
theory, it is similar to the Poincaré-Birkhoff local normal forms for vector fields
(Levi normal forms are governed by semi-simple group actions while Poincaré-Birkhoff
normal forms are governed by torus actions, see
\cite{Zung-Birkhoff2001,Zung-PD2002}.

5. Theorem \ref{thm:LeviPoisson} provides an useful tool for the study of
linearization of Poisson structures. Using it, J.-P. Dufour and I recently showed in
\cite{DuZu-affn2002} that the Lie algebra $\mathfrak{aff}(n, \bbK)$ of infinitesimal
affine transformations of $\bbK^n$ is analytically nondegenerate.
\\

It is natural that not only Poisson structures but other geometric structures
related to infinite-dimensional Lie algebras admit formal or analytic Levi
decomposition as well. For example, Cerveau \cite{Cerveau} showed the existence of a
formal Levi decomposition for singular foliations.\footnote{As far as we know, the
existence of an analytic Levi decomposition for singular foliations remains an open
problem.} In this paper, we will show that analytic Lie algebroids also admit local
analytic Levi decomposition.

Recall (see e.g.
\cite{CaWe,Dufour-algebroid,Fernandes-Algebroid2001,Weinstein-algebroid}) that a
smooth Lie algebroid over a manifold $M$ is a vector bundle $A\rightarrow M$ with a
Lie algebra structure on its space $\Gamma (A)$ of smooth sections  and a bundle map
$\# : A \rightarrow TM$ (called the anchor) inducing a Lie algebra homomorphism from
sections of $A$ to vector fields on $M$, such that $[s,fs']=f[s,s']+ (\# s \cdot f)
s'$ for sections $s$ and $s'$ and functions $f$. In the analytic category, one
replaces $\Gamma(A)$ by the sheaf of local analytic sections. A point $x \in M$ is
called singular for the algebroid $A$ if the rank of the anchor map $\#_x : A_x
\rightarrow T_xM$ (where $A_x$ is the fiber of $A$ over $x$) is smaller than at
other points. Due to the local splitting theorem for Lie algebroids (see
\cite{Dufour-algebroid,Fernandes-Algebroid2001,Weinstein-algebroid}), in the study
of local normal forms of Lie algebroids near a singular point $x$, we may assume
that the rank of $\#_x : A_x \rightarrow T_xM$ is zero.

Let $A$ be a local analytic Lie algebroid of dimension $N$ over $(\bbK^n,0)$ such
that the anchor map $\# : A_x \rightarrow T_x\bbK^n$ vanishes at $x= 0$. Denote by
$s_1,...,s_N$ an analytic local basis of sections of $A$, and $(x_1,...,x_n)$ an
analytic local system of coordinates of $(\bbK^n,0)$. Then we have $[s_i,s_j] =
\sum_k c_{ij}^k s_k +$ higher order terms in $s_1,...,s_N$, and $\# s_i = \sum_{j,k}
b_{ij}^k x_k \partial / \partial x_j +$ higher order terms in $x_1,...,x_n$. If we
forget about the terms of order greater or equal to 2, then we get an
$N$-dimensional Lie algebra with structural coefficients $c_{ij}^k$, which acts on
$\bbK^n$ via linear vector fields $\sum_{j,k} b_{ij}^k x_k \partial /
\partial x_j$. (The action Lie algebroid associated to this linear Lie algebra action
is called the linear part of the algebroid $A$ at $0$). Denote this $N$-dimensional
Lie algebra by $\frak L$, and by $\frak L = \fg \ltimes \frak r$ its Levi
decomposition. We are looking for a Levi factor of $\Gamma(A)$, where $\Gamma(A)$
now denotes the infinite-dimensional Lie algebra of local analytic sections of $A$
(the Lie bracket is given by the algebroid structure of $A$), i.e. a subalgebra of
$\Gamma(A)$ which is isomorphic to $\fg$. Once such a Levi factor is found, its
action on the algebroid $A$ can be linearized by Guillemin-Sternberg theorem
\cite{GuSt-Linearization1968}, because $\fg$ is semi-simple.

\begin{thm}
\label{thm:LeviAlgebroid} Let $A$ be a local $N$-dimensional analytic Lie algebroid
over $(\bbK^n,0)$ with the anchor map $\# : A \to T\bbK^n$, such that $\#a = 0$ for
any $a \in A_0$, the fiber of $A$ over point $0$. Denote by $\frak L$ the
$N$-dimensional Lie algebra in the linear part of $A$ at $0$, and by $\frak L = \fg
\ltimes \frak r$ its Levi decomposition. Then there exists a local analytic system
of coordinates $(x^\infty_1,...,x^\infty_n)$ of $(\bbK^n,0)$, and a local analytic
basis of sections
$(s^\infty_1,s^\infty_2,...,s^\infty_m,v^\infty_1,...,v^\infty_{N-m})$ of $A$, where
$m = \dim \fg$, such that we have :
\begin{equation}
\begin{array}{l}
[s^\infty_i,s^\infty_j] = \sum_k c_{ij}^k s^\infty_k  \ , \cr [s^\infty_i,
v^\infty_j] = \sum_k a_{ij}^k v^\infty_k \ ,\cr \# s^\infty_i = \sum_{j,k} b_{ij}^k
x^\infty_k \partial /
\partial x^\infty_j \ ,
\end{array}
\end{equation}
where $c_{ij}^k,a_{ij}^k,b_{ij}^k$ are constants, with $c_{ij}^k$ being the
structural coefficients of the semi-simple Lie algebra $\fg$.
\end{thm}

{\it Remarks.}

1. In the above theorem, when $\frak L = \fg$, we get the analytic linearization of
Lie algebroids with semi-simple linear part. The formal version of this
linearization result has been obtained by Dufour \cite{Dufour-algebroid} and
Weinstein \cite{Weinstein-algebroid}. When the Lie algebroid is an action algebroid,
we also recover classical results about the linearization of analytic actions of
semi-simple Lie groups and Lie algebras.

2. The proof of the above theorem is absolutely similar to that of Theorem
\ref{thm:LeviPoisson}. In fact, since Lie algebroid structures on a vector bundle
may be viewed as ``fiberwise linear'' Poisson structures on the dual bundle (see
e.g. \cite{CaWe}), Theorem \ref{thm:LeviAlgebroid} may be viewed as a special case
of Theorem \ref{thm:LeviPoisson}.
\\

The rest of this paper is devoted to the proof of Theorem \ref{thm:LeviPoisson} and
Theorem \ref{thm:LeviAlgebroid}. We will first prove Theorem \ref{thm:LeviPoisson},
and then show a few modifications to be made to our proof of Theorem
\ref{thm:LeviPoisson} to get a proof of Theorem \ref{thm:LeviAlgebroid}.

\section{Formal Levi decomposition}
\label{section:formal}

In this section we will construct by recurrence a formal system of coordinates
$(x^{\infty}_1,...,x^{\infty}_m,y^{\infty}_1,...,y^{\infty}_{n-m})$ which satisfy
Relations (\ref{eqn:LeviPoisson2}) for a given local analytic Poisson structure
$\Pi$. We will later use analytic estimates to show that our construction actually
yields a local analytic system of coordinates. Let us mention that our construction
below of the Levi decomposition differs from the constructions of Wade \cite{Wade} and
Weinstein \cite{Weinstein-algebroid}. Their constructions are simpler and are good
enough to show the existence of a formal Levi factor. However, in order to show the
existence of an analytic Levi factor (using the fast convergence method), we need a
more complicated construction, in which each step in a recurrence process consists
of 2 substeps: the first substep is to find an ``almost Levi factor''. The second
substep consists of ``almost linearizing'' this ``almost Levi factor''.

We begin the first step with the original linear coordinate system
$$(x^{0}_1,...,x^{0}_m,y^{0}_1,...,y^{0}_{n-m}) = (x_1,...,x_m,y_1,...,y_{m-n}).$$
For each positive integer $l$, after Step $l$ we will find a local coordinate system
$(x^{l}_1,...,x^{l}_m,y^{l}_1,...,y^{l}_{n-m})$ with the following properties
(\ref{eqn:xil}), (\ref{eqn:phil}), (\ref{eqn:Xil}) :

\begin{equation}
\label{eqn:xil} (x^{l}_1,...,x^{l}_m,y^{l}_1,...,y^{l}_{n-m}) =
(x^{l-1}_1,...,x^{l-1}_m,y^{l-1}_1,...,y^{l-1}_{n-m}) \circ \phi_{l} \ ,
\end{equation}
where $\phi_l$ is a local analytic diffeomorphism of $(\bbK^n,0)$ of the type
\begin{equation}
\label{eqn:phil} \phi_l(x) = x + \psi_l (x) \ , \ \psi_l(x) \in O(|x|^{2^{l-1}+1})
\end{equation}
(i.e., $\psi_l(x)$ contains only terms of order greater or equal to $2^{l-1}+1$).

Denote by
\begin{equation}
X^{l}_i = X_{x^{l}_i} \ (i = 1,...,m)
\end{equation}
the Hamiltonian vector field of $x^{l}_i$ with respect to our Poisson structure
$\Pi$. Then we have
\begin{equation}
X^l_i = \hat{X}^l_i + Y^l_i \ ,
\end{equation}
where
\begin{equation}
\label{eqn:Xil} \hat{X}^{l}_i = \sum_{jk} c_{ij}^k x^{l}_k \frac{\partial }{
\partial x^{l}_j} + \sum_{jk} a_{ij}^k y^{l}_k \frac{\partial }{
\partial y^{l}_j} \ , \ Y^{l}_i \in O(|x|^{2^l+1}) \ ,
\end{equation}
i.e., $\hat{X}^l_i$ is the linear part of $X^{l}_i = X_{x^{l}_i}$ in the coordinate
system $(x^l_1,...,y_{n-m}^l)$, $c_{ij}^k$ and $a_{ij}^k$ are structural constants
as appeared in Theorem \ref{thm:LeviPoisson}, and $Y^l_i = X^l_i - \hat{X}^l_i$ does
not contain terms of order $\leq 2^l$.

Of course, when $l=0$, then Relation (\ref{eqn:Xil}) is satisfied by the assumptions
of Theorem \ref{thm:LeviPoisson}. Let us show how to construct the coordinate system
$(x^{l+1}_1,...,y^{l+1}_{n-m})$ from the coordinate system $(x^l_1,...,y_{n-m}^l)$.
Denote

\begin{equation}
\cO_l = \{ {\rm local \ analytic \ functions \ in \ (\bbK^n,0) \ without \ terms \
of \ oder} \leq 2^l \} \ .
\end{equation}

Due to Relations (\ref{eqn:xil}) and (\ref{eqn:phil}), it doesn't matter if we use
the coordinate system $(x_1,...,x_m,y_1,...,y_{n-m})$ or the coordinate system
$(x^{l}_1,...,x^{l}_m,y^{l}_1,...,y^{l}_{n-m})$ in the above definition of $\cO_l$.
It follows from Relation (\ref{eqn:Xil}) that

\begin{equation}
\label{eqn:fij} f^{l}_{ij} := \{x^{l}_i, x^{l}_j \} - \sum_k c_{ij}^k x^{l}_k =
Y^l_i(x^l_j) \in \cO_l \ .
\end{equation}

Denote by $(\xi_1,...,\xi_m)$ a fixed basis of the semi-simple algebra $\fg$, with

\begin{equation}
[\xi_i,\xi_j] = \sum_k c_{ij}^k \xi_k \ .
\end{equation}

Then $\fg$ acts on $\cO$ via vector fields $\hat{X}^{l}_1,...,\hat{X}^{l}_m$, and
this action induces the following linear action of $\fg$ on the finite-dimensional
vector space $\cO_l/\cO_{l+1}$ : if $g \in \cO_l$, considered modulo $\cO_{l+1}$,
then we put

\begin{equation}
\label{eqn:action} \xi_i \cdot g := \hat{X}^{l}_i (g) = \sum_{jk} c_{ij}^k x^l_k
\frac{\partial g}{\partial x^l_j} +  \sum_{jk} a_{ij}^k y^l_k \frac{\partial
g}{\partial y^l_j} \ \  mod \ \cO_{l+1} \ .
\end{equation}

Notice that if $g \in \cO_l$ then $Y^l_i(g) \in \cO_{l+1}$, and hence we have
\begin{equation}
\xi_i \cdot g = X^l(g) \ mod \ \cO_{l+1} \ = \{x^l_i, g\} \ mod \ \cO_{l+1} \ .
\end{equation}

The functions $f^{l}_{ij}$ in (\ref{eqn:fij}) form a 2-cochain $f^{l}$ of $\fg$ with
values in the $\fg$-module $\cO_l/\cO_{l+1}$ :
\begin{equation}
\begin{array}{c}
f^{l} : \fg \wedge \fg \to \cO_l / \cO_{l+1} \\
f^{l}(\xi_i \wedge \xi_j) := f^{l}_{ij} \ mod \ \cO_{l+1} \ = \{x^{l}_i, x^{l}_j \}
- \sum_k c_{ij}^k x^{l}_k \ mod \ \cO_{l+1} \ .
\end{array}
\end{equation}
In other words, if we denote by $\fg^{\ast}$ the dual space of $\fg$, and by
$(\xi_1^{\ast},...,\xi_m^{\ast})$ the basis of $\fg^{\ast}$ dual to
$(\xi_1,...,\xi_m)$, then we have
\begin{equation}
f^l = \sum_{i < j} \xi^{\ast}_i \wedge \xi^{\ast}_j \otimes (f^l_{ij} \ mod \
\cO_{l+1}) \in \wedge^2 \fg^{\ast} \otimes \cO_l/\cO_{l+1} \ .
\end{equation}

It follows from (\ref{eqn:fij}), and the Jacobi identity for the Poisson bracket of
$\Pi$ and the algebra $\fg$, that the above 2-cochain is a 2-cocycle. Because $\fg$
is semi-simple, we have $H^2 (\fg, \cO_l/\cO_{l+1}) = 0$, i.e. the second cohomology
of $\fg$ with coefficients in $\fg$-module $\cO_l/\cO_{l+1}$ vanishes, and therefore
the above 2-cocycle is a coboundary. In other words, there is an 1-cochain
\begin{equation}
w^{l} \in \fg^{\ast} \otimes \cO_l/\cO_{l+1}
\end{equation}
such that
\begin{equation}
\label{eqn:w} f^{l}(\xi_i \wedge \xi_j) = \xi_i \cdot w^{l}(\xi_j) - \xi_j \cdot
w^{l}(\xi_i) - w^{l}(\sum_k c_{ij}^k\xi_k) \ .
\end{equation}

Denote by $w^{l}_i$ the element of $\cO_l$ which is a polynomial of order $\leq
2^{l+1}$ in variables $(x^{l}_1,...,x^{l}_m,y^{l}_1,y^{l}_{n-m})$ such that the
projection of $w^{l}_i$ in $\cO_l/\cO_{l+1}$ is $w^{l}(\xi_i)$. Define $x^{l+1}_i$
as follows:

\begin{equation}
x^{l+1}_i = x^{l}_i - w^{l}_i \ \ (i = 1, \dots, m) \ .
\end{equation}

Then it follows from (\ref{eqn:fij}) and (\ref{eqn:w}) that we have

\begin{equation}
\label{eqn:wi2} \{x^{l+1}_i, x^{l+1}_j\} - \sum_k c_{ij}^k (x^{l+1}_k) \in \cO_{l+1}
\ for \ i,j \leq m \ .
\end{equation}

Denote by $\cY^l$ the space of local analytic vector fields of the type $u =
\sum_{i=1}^{n-m} u_i \partial / \partial y^l_i$ (with respect to the coordinate
system $(x^l_1,...,y_{n-l}^l)$), with $u_i$ being local analytic functions. For each
natural number $k$, denote by $\cY^l_k$ the following subspace of $\cY^l$:

\begin{equation}
\cY^l_k = \left\{ \Big. u = \sum_{i=1}^{n-m} u_i \partial / \partial y^l_i \ \ \Big|
\ \ u_i \in \cO_k \right\} \ .
\end{equation}

Then $\cY^l$, as well as $\cY^l_l/\cY^l_{l+1}$, are $\fg$-modules under the
following action :

\begin{equation}
\xi_i \cdot \sum_j u_j \partial / \partial y^l_j := [\hat{X}^l_i, u] = \Big[
\sum_{jk} c_{ij}^k x^l_k \frac{\partial }{
\partial x^l_j} + \sum_{jk} a_{ij}^k y^l_k \frac{\partial }{
\partial y^l_j}\ ,\ \sum_j u_j \partial /
\partial y^l_j \Big] \ .
\end{equation}

The above linear action of $\fg$ on $\cY_l / \cY_{l+1}$ can also be written as
follows :
\begin{equation}
\xi_i \cdot \sum_j u_j \partial / \partial y^l_j = \sum_{j} (\{x_i^l,u_j\} - \sum_k
a_{ij}^k u_k) \partial / \partial y^l_j \ \ mod \ \ \cY^l_{l+1} \ .
\end{equation}

Define the following 1-cochain of $\fg$ with values in $\cY^l_l/\cY^l_{l+1}$ :
\begin{equation}
\label{eqn:1cocycle} \sum_{i=1}^m  \big( \xi_i^{\ast} \otimes \big( \sum_{j=1}^{n-m}
(\{x^{l+1}_i,y^l_j\} - \sum_k a_{ij}^k y^l_k)
\partial / \partial y^l_j \ mod \ \cY^l_{l+1} \big) \big) \ \in \ \fg^{\ast} \otimes
\cY^l_l/\cY^l_{l+1} \ .
\end{equation}

Due to Relation (\ref{eqn:wi2}), the above 1-cochain is an 1-cocycle. Since $\fg$ is
semi-simple, we have $H^1(\fg,\cY^l_l / \cY^l_{l+1}) = 0$, and the above 1-cocycle
is an 1-coboundary. In other words, there exists a vector field $\sum_{j=1}^{n-m}
v^l_j
\partial / \partial y^l_j \in \cY^l_l$, with $v^l_j$ being a polynomial function of degree
$\leq 2^{l+1}$ in variables $(x^l_1,...,y^l_{n-m})$, such that for every $i=
1,...,m$ we have
\begin{equation}
\label{eqn:vi} \sum_j (\{x^{l+1}_i,y^l_j\} - \sum a_{ij}^k y^l_k) \partial/
\partial y^l_j =
\sum_j ( \{x^l_i,v^l_j\} - \sum a_{ij}^k v^l_k) \partial/
\partial y^l_j \ \  mod \ \ \cY^l_{l+1} \ .
\end{equation}

We now define the new system of coordinates as follows :

\begin{equation}
\begin{array}{l}
x^{l+1}_i = x^l_i - w^l_i \ (i=1,...,m), \cr y^{l+1}_i = y^{l}_i - v^l_i \
(i=1,...,n-m),
\end{array}
\end{equation}
where functions $w^l_i,v^l_i \in \cO_l$ are chosen as above. In particular,
Relations (\ref{eqn:wi2}) and (\ref{eqn:vi}) are satisfied, which means that

\begin{equation}
\begin{array}{l}
\{x_i^{l+1},x_j^{l+1}\} - \sum c_{ij}^k x_k^{l+1} \in \cO_{l+1} \ , \\
\{x_i^{l+1},y_j^{l+1}\} - \sum a_{ij}^k y_k^{l+1} \in \cO_{l+1} \ ,
\end{array}
\end{equation}
i.e. Relation (\ref{eqn:Xil}) is satisfied with $l$ replaced by $l+1$. Of course,
Relations (\ref{eqn:xil}) and (\ref{eqn:phil}) are also satisfied with $l$ replaced
$l+1$, with $\phi_{l+1} = Id + \psi_{l+1}$ and
\begin{equation}
\psi_{l+1} = - (w^l_1,...,w^l_m,v^l_1,...,v^l_{n-m}) \in (O_l)^{n} \ .
\end{equation}
Recall that, by the above construction, the functions
$w^l_1,...,w^l_m,v^l_1,...,v^l_{n-m}$ are polynomial functions of degree $\leq
2^{l+1}$ in variables $(x^l_1,...,y^l_{n-m})$, which do not contain terms of degree
$\leq 2^l$.

Define the following limits
\begin{equation}
\begin{array}{l}
(x^{\infty}_1,...,y^{\infty}_{n-m}) = \lim_{l \to \infty} (x^{l}_1,...,y^{l}_{n-m})
\ , \cr \Phi_{\infty} = \lim_{l \to \infty} \Phi_l \ \ {\rm where} \ \Phi_l = \phi_1
\circ ... \circ \phi_l \ .
\end{array}
\end{equation}

It is clear that the above limits exist in the formal category,
$(x^{\infty}_1,...,y^{\infty}_{n-m}) = (x^{0}_1,...,y^{0}_{n-m}) \circ
\Phi_{\infty}$, and the formal coordinate system
$(x^{\infty}_1,...,y^{\infty}_{n-m})$ satisfies Relation (\ref{eqn:LeviPoisson2}).
To prove Theorem \ref{thm:LeviPoisson}, it remains to show that we can choose
functions $w_i^l,v_i^l$ in such a way that $(x^{\infty}_1,...,y^{\infty}_{n-m})$ is
in fact a local analytic system of coordinates.

\section{Normed vanishing of cohomologies}

In this section, using ``normed vanishing'' of first and second cohomologies of
$\fg$, we will obtain some estimates on $w^l_i = x^l_i - x^{l+1}_i$ and $v^l_i =
y^l_i - y^{l+1}_i$. See e.g. \cite{Jacobson} for some basic results on semi-simple
Lie algebras and their representations which will be used below.

We will denote by $\fg_{\bbC}$ the algebra $\fg$ if $\bbK = \bbC$, and the
complexification of $\fg$ if $\bbK = \bbR$. So $\fg_\bbC$ is a complex semi-simple
Lie algebra of dimension $m$. Denote by $\fg_0$ the compact real form of
$\fg_{\bbC}$, and identify $\fg_{\bbC}$ with $\fg_0 \otimes_\bbR \bbC$. Fix an
orthonormal basis $(e_1,...,e_m)$ of $\fg_{\bbC}$ with respect to the Killing form :
$<e_i,e_j> = \delta_{ij}$. We may assume that $e_1,...,e_m \in \sqrt{-1} \fg_0$.
Denote by $\Gamma = \sum_i e_i^2$ the Casimir element of $\fg_{\bbC}$ : $\Gamma$
lies in the center of the universal enveloping algebra ${\mathcal U}(\fg_{\bbC})$
and does not depend on the choice of the basis $(e_i)$. When $\bbK = \bbR$ then
$\Gamma$ is real, i.e., $\Gamma \in {\mathcal U}(\fg)$.

Let $W$ be a finite dimensional complex linear space endowed with a Hermitian metric
denoted by $<,>$. If $v \in W$ then its norm is denoted by $\|v\| = \sqrt{<v,v>}$.
Assume that $W$ is a Hermitian $\fg_0$-module. In other words, the linear action of
$\fg_0$ on $W$ is via infinitesimal unitary (i.e. skew-adjoint) operators. $W$ is a
$\fg_{\bbC}$-module via the identification $\fg_{\bbC} = \fg_0 \otimes_\bbR \bbC$.
We have the decomposition $W = W_0 + W_1$, where $W_1 = \fg_{\bbC} \cdot W$ (the
image of the representation), and $\fg_{\bbC}$ acts trivially on $W_0$. Since $W_1$
is a $\fg_{\bbC}$-module, it is also a ${\mathcal U}(\fg_{\bbC})$-module. The action
of $\Gamma$ on $W_1$ is invertible : $\Gamma \cdot W_1 = W_1$, and we will denote by
$\Gamma^{-1}$ the inverse mapping.

Denote by $\fg^{\ast}_\bbC$ the dual of $\fg_\bbC$, and by
$(e^{\ast}_1,...,e^{\ast}_m)$ the basis of $\fg^{\ast}_\bbC$ dual to
$(e_1,...,e_m)$. If $w \in \fg^{\ast}_\bbC \otimes W$ is an 1-cochain and $f :
\wedge^2 \fg^{\ast}_\bbC \otimes W$ is a 2-cochain with values in $W$, then we will
define the norm of $f$ and $w$ as follows :

\begin{equation}
\|w\| = \max_i \|w(e_i)\| \ , \ \|f\| = \max_{i,j} \|f(e_i \wedge e_j)\| \ .
\end{equation}

Since $H^2(\fg, \bbK) = 0$, there is a (unique) linear map $h_0 :  \wedge^2 \fg^\ast
\rightarrow \fg^\ast$ such that if $u \in \wedge^2 \fg^\ast$ is a 2-cocycle for the
trivial representation of $\fg$ in $\bbK$ (i.e. $u([x,y],z) + u([y,z],x) +
u([z,x],y) = 0$ for any $x,y,z \in \fg$), then $u = \delta h_0(u)$, i.e. $u(x,y) =
h_0(u) ([x,y])$. By complexifying $h_0$ if $\bbK = \bbR$, and taking its tensor
product with the projection map $P_0 : W \to W_0$, we get a map

\begin{equation}
h_0 \otimes P_0 : \wedge^2 \fg^\ast_\bbC \otimes W \rightarrow  \fg^\ast_\bbC
\otimes W_0 \ .
\end{equation}

Define another map
\begin{equation}
h_1 : \wedge^2 \fg^\ast_\bbC \otimes W \rightarrow  \fg^\ast_\bbC \otimes W_1
\end{equation}
as follows : if $f \in \wedge^2 \fg^{\ast}_\bbC \otimes W$ then we put
\begin{equation}
\label{eqn:homotopy1} h_1 (f)  =  \sum_i e_i^\ast \otimes (\Gamma^{-1} \cdot \sum_j
(e_j \cdot f(e_i \wedge e_j ))) \ .
\end{equation}
Then the map
\begin{equation}
h = h_0 \otimes P_0 + h_1 : \wedge^2 \fg^\ast_\bbC \otimes W \rightarrow
\fg^\ast_\bbC \otimes W
\end{equation}
is an explicit homotopy operator, in the sense that if $f \in \wedge^2
\fg^{\ast}_\bbC \otimes W$ is a 2-cocycle (i.e. $\delta f = 0$ where $\delta$
denotes the differential of the Eilenberg-Chevalley complex $ ... \rightarrow
\wedge^k \fg^{\ast}_\bbC \otimes W \rightarrow \wedge^{k+1} \fg^{\ast}_\bbC \otimes
W \rightarrow ... $), then $f = \delta (h(f))$.

Similarly, the map $h : \fg^\ast_\bbC \otimes W \rightarrow W$ defined by
\begin{equation}
\label{eqn:homotopy2}
h (w) = \Gamma^{-1} \cdot (\sum_i e_i \cdot w(e_i))
\end{equation}
is also a homotopy operator, in the sense that if $w \in \fg^\ast_\bbC \otimes W$ is
an 1-cocycle then $w = \delta (h(w))$.

When $\bbK = \bbR$, i.e. when $\fg_\bbC$ is the complexification of $\fg$, then the
above homotopy operators $h$ are real, i.e. they map real cocycles into real
cochains.

The above formulas make it possible to control the norm of a primitive of a
1-cocycle $w$ or a 2-cocycle $f$ in terms of the norm of $w$ or $f$ : we have the
following lemma, which has been (essentially) proved by Conn in Proposition 2.1 of
ref. \cite{Conn} and Proposition 2.1 of ref. \cite{Conn-Smooth1985}.

\begin{lemma}
\label{lemma:inverse} There is a positive constant $D$ (which depends on $\fg$ but
does not depend on $W$) such that with the above notations we have
\begin{equation}
\|h(f)\| \leq D \|f\| \ {\rm and} \ \|h(w)\| \leq D \|w\|
\end{equation}
for any 1-cocycle $w$ and any 2-cocycle $f$ of $\fg_\bbC$ with values in $W$.
\end{lemma}

{\it Proof}. (See Proposition 2.1 of \cite{Conn} and Proposition 2.1 of
\cite{Conn-Smooth1985}) We can decompose $W$ into an orthogonal sum (with respect to
the Hermitian metric of $W$) of irreducible modules of $\fg_0$. The above homotopy
operators decompose correspondingly, so it is enough to prove the above lemma for
the case when $W$ is a non-trivial irreducible module, which we will now suppose.
Let $\lambda \neq 0$ denote the highest weight of the irreducible $\fg_0$-module
$W$, and by $\delta$ one-half the sum of positive roots of $\fg_0$ (with respect to
a fixed Cartan subalgebra and Weil chamber). Then $\Gamma$ acts on $W$ by
multiplication by the scalar $\langle\lambda, \lambda + 2\delta\rangle$, which is
greater or equal to $\|\lambda\|^2$. Denote by $\mathcal J$ the weight lattice of
$\fg_0$, and $D = m ( \min_{\gamma \in \mathcal J} \|\gamma\|)^{-1}$. Then $D <
\infty$ does not depend on $W$, and $\|\lambda\|^2 > \frac{m \|\lambda\|}{D}$, which
implies that the norm of the inverse of the action of $\Gamma$ on $W$ is smaller or
equal to $\frac{D}{m \|\lambda\|}$. On the other hand, the norm of the action of
$e_i$ on $W$ is smaller or equal to $\|\lambda\|$ for each $i = 1,...,m$ (recall
that $\sqrt{-1}e_i \in \fg_0$ and $<e_i,e_i> = 1$), hence the norm of the operator
$\sum_{i=1}^m e_i \cdot \Gamma^{-1} : W \rightarrow W$ is smaller or equal to $D$.
Now apply Formulas (\ref{eqn:homotopy1}) and (\ref{eqn:homotopy2}). The
lemma is proved. \hfill $\square$ \\

Let us now apply the above lemma to $\fg$-modules $\cO_l/\cO_{l+1}$ and
$\cY^l_l/\cY^l_{l+1}$ introduced in the previous section. Recall that $\fg$ is a
Levi factor of $\frak L$, the space of linear functions in $\bbK^n$, which is a Lie
algebra under the linear Poisson bracket $\Pi_1$. $\fg$ acts on $\frak L$ by the
(restriction of the) adjoint action, and on $\bbK^n$ by the coadjoint action. By
complexifying these actions if necessary, we get a natural action of $\fg_\bbC$ on
$(\bbC^n)^{\ast}$ (the dual space of $\bbC^n$) and on $\bbC^n$. The elements
$x_1,...,x_m,y_1,...,y_{n-m}$ of the original linear coordinate system in $\bbK^n$
may be view as a basis of $(\bbC^n)^{\ast}$. Notice that the action of $\fg_\bbC$ on
$(\bbC^n)^{\ast}$ preserves the subspace spanned by $(x_1,...,x_m)$ and the subspace
spanned by $(y_1,...,y_{n-m})$. Fix a basis $(z_1,...,z_n)$ of $(\bbC^n)^{\ast}$,
such that the Hermitian metric of $(\bbC^n)^{\ast}$ for which this basis is
orthonormal is preserved by the action of $\fg_0$, and such that
\begin{equation}
z_{i} = \sum_{j \leq m} A_{ij} x_j + \sum_{j \leq n-m} A_{i,j+m} y_j \ ,
\end{equation}
with the constant transformation matrix $(A_{ij})$ satisfying the following
condition :
\begin{equation}
\label{eqn:A} A_{ij} = 0 \ {\rm if} \ (i \leq m < j \ {\rm or} \ j \leq m < i) \ .
\end{equation}
Such a basis $(z_1,...,z_n)$ always exists, and we may view $(z_1,...,z_n)$ as a
linear coordinate system on $\bbC^n$. We will also define local complex analytic
coordinate systems $(z^l_1,...,z^l_n)$ as follows :

\begin{equation}
\label{eqn:zl} z^l_{i} = \sum_{j \leq m} A_{ij} x^l_j + \sum_{j \leq n-m} A_{i,j+m}
y^l_j  \ .
\end{equation}

Let $l$ be a natural number, $\rho$ a positive number, and $f$ a local complex
analytic function of $n$ variables. Define the following ball $B_{l,\rho}$ and
$L^2$-norm $\|f\|_{l,\rho}$, whenever it makes sense :

\begin{equation}
\label{eqn:ball} B_{l,\rho} = \Big\{ x \in \bbC^n \ | \ \sqrt{\sum |z^l_i(x)|^2}
\leq \rho \Big\} \ ,
\end{equation}

\begin{equation}
\|f\|_{l,\rho} = \sqrt{\frac{1}{V_\rho}\int_{B_{l,\rho}} |f(x)|^2 d\mu_l} \ ,
\end{equation}
where $d\mu_l$ is the standard volume form in the complex ball $B_{l,\rho}$ with
respect to the coordinate system $(z^l_1,...,z^l_n)$, and $V_\rho$ is the volume of
$B_{l,\rho}$, i.e. of an $n$-dimensional complex ball of radius $\rho$.

We will say that the ball $B_{l,\rho}$ is well-defined if it is analytically
diffeomorphic to the standard ball of radius $\rho$ via the coordinate system
$(z^l_1,...,z^l_n)$, and will use $\|f\|_{l,\rho}$ only when $B_{l,\rho}$ is
well-defined. When $B_{l,\rho}$ is not well-defined we simply put $\|f\|_{l,\rho} =
\infty $. We will write $B_\rho$ and $\|f\|_{\rho}$ for $B_{0,\rho}$ and
$\|f\|_{0,\rho}$ respectively. If $f$ is a real analytic function (the case when
$\bbK =\bbR$), we will complexify it before taking the norms.

It is well-known that the $L^2$-norm $\|f\|_{\rho}$ is given by a Hermitian metric,
in which the monomial functions form an orthogonal basis : if $f = \sum_{\alpha \in
{\mathbb N}^n } a_\alpha \prod_i z_i^{\alpha_i}$ and $g = \sum_{\alpha \in {\mathbb
N}^n } b_\alpha \prod_i z_i^{\alpha_i}$ then the scalar product $\langle f,g
\rangle_\rho$ is given by
\begin{equation}
\label{eqn:product} \langle f,g \rangle_\rho = \sum_{\alpha \in {\mathbb N}^n}
\frac{\alpha!(n-1)!}{(|\alpha| + n - 1 )!}\rho^{2|\alpha|} a_\alpha \bar{b}_\alpha \
,
\end{equation}
(where $\alpha! = \prod_i \alpha_i !, |a| = \sum \alpha_i$, and $\bar{b}$ is the
complex conjugate of $b$), and the norm $\|f\|_{\rho}$ is given by
\begin{equation}
\|f\|_{\rho} = \left( \sum_{\alpha \in {\mathbb N}^n} \frac{\alpha!(n-1)!}{(|\alpha|
+ n - 1 )!}|c_\alpha|^2\rho^{2|\alpha|} \right)^{1/2} \ .
\end{equation}

The above scalar product turns $\cO_l/\cO_{l+1}$ into a Hermitian space, if we
consider elements of $\cO_l/\cO_{l+1}$ as polynomial functions of degree less or
equal to $2^{l+1}$ and which do not contain terms of order $\leq 2^l$. Of course,
when $\bbK = \bbR$ we will have to complexify $\cO_l/\cO_{l+1}$, but will redenote
$(\cO_l/\cO_{l+1})_\bbC$ by $\cO_l/\cO_{l+1}$, for simplicity.

Similarly, for the space $\cY^l$ of local vector fields of the type $u =
\sum_{i=1}^{n-m} u_i \partial / \partial z^l_{i+m}$ (due to (\ref{eqn:A}) and
(\ref{eqn:zl}), this is the same as the space of vector fields of the type
$\sum_{i=1}^{n-m} u'_i
\partial /
\partial y^l_{i}$ defined in the previous section, up to a complexification if
$\bbK = \bbR$), we define the $L^2$-norms as follows :
\begin{equation}
\|u\|_{l,\rho} = \sqrt{\frac{1}{V_\rho}\int_{B_{l,\rho}} \sum_{i=1}^{n-m} |u_i(x)|^2
d\mu_l} \ .
\end{equation}
These $L^2$-norms are given by Hermitian metrics similar to (\ref{eqn:product}),
which make $\cY^l_l/\cY^l_{l+1}$ into Hermitian spaces.

Remark that if $u = (u_1,...,u_{n-m})$ then
\begin{equation}
\sum_i \|u_i\|_{l,\rho} \geq \|u\|_{l,\rho} \geq \max_i \|u_i\|_{l,\rho} \ .
\end{equation}

It is an important observation that, since the action of $\fg_{0}$ on $\bbC^n$
preserves the Hermitian metric of $\bbC^n$, its actions on $\cO_l/\cO_{l+1}$ and
$\cY^l_l/\cY^l_{l+1}$, as given in the previous section, also preserve the Hermitian
metrics corresponding to the norms $\|f\|_{l,\rho}$ and $\|u\|_{l,\rho}$ (with the
same $l$). Thus, applying Lemma \ref{lemma:inverse} to these $\fg_\bbC$-modules, we
get :

\begin{lemma}
\label{lemma:inverse2}
 There is a positive constant $D_1$ such that for any $l \in \bbN$ and
any positive number  $\rho$ there exist local analytic functions
$w^l_1,...,w^l_m,v^l_1,...,v^l_{n-m}$, which satisfy the relations of the previous
section, and which have the following additional property whenever $B_{l,\rho}$ is
well-defined :
\begin{equation}
\max_i \|w^l_i\|_{l,\rho} \leq D_1 . \max_{i,j} \|\{x^l_i,x^l_j\} - \sum_k c_{ij}^k
x^l_k \|_{l,\rho}
\end{equation}
and
\begin{equation}
\max_i \|v^l_i\|_{l,\rho} \leq D_1 . \max_{i,j} \|\{x^l_i-w^l_i,y^l_j\} - \sum_k
a_{ij}^k y^l_k \|_{l,\rho} \ .
\end{equation}
\end{lemma} \hfill $\square$

\section{Proof of convergence}

Besides the $L^2$-norms defined in the previous section, we will need the following
$L^{\infty}$-norms : If $f$ is a local function then put
\begin{equation}
|f|_{l,\rho} = \sup_{x \in B_{l,\rho}} |f(x)| \ ,
\end{equation}
where the complex ball $B_{l,\rho}$ is defined by (\ref{eqn:ball}). Similarly, if $g
= (g_1,...,g_N)$ is a vector-valued local map then put $|g|_{l,\rho} = \sup_{x \in
B_{l,\rho}} \sqrt{\sum_i |g_i(x)|^2}$. For simplicity, we will write $|f|_\rho$ for
$|f|_{0,\rho}$.

For the Poisson structure $\Pi$, we will use the following norms :
\begin{equation}
|\Pi|_{l,\rho} := \max_{i,j=1,...,n} \{ |\{z^l_i,z^l_j\}|_{l,\rho} \} \ .
\end{equation}

Due to the following lemma, we will be able to use the norms $|f|_{\rho}$ and
$\|f\|_{\rho}$ interchangeably for our purposes, and control the norms of the
derivatives :

\begin{lemma}
\label{lemma:2norms} For any $\epsilon > 0$ there is a finite number $K < \infty$
depending on $\epsilon$ such that for any integer $l > K$, positive number $\rho$,
and local analytic function $f \in \cO_l$ we have
\begin{equation}
\label{eqn:changeradius} |f|_{(1 + \epsilon/l^2)\rho} \geq exp(2^{l/2}) |f|_{(1 +
\epsilon/2l^2)\rho} \geq \rho |df|_\rho \ ,
\end{equation}
and
\begin{equation}
\label{eqn:2norms} |f|_{(1 - \epsilon/l^2)\rho} \leq \|f\|_{\rho} \leq |f|_{\rho} \
.
\end{equation}
\end{lemma}

The above lemma, and other lemmas in this section, will be proved in the subsequent
section.

The key point in the proof of Theorem \ref{thm:LeviPoisson} is the following
proposition.

\begin{prop}
\label{prop:mainPoisson} Under the assumptions of Theorem \ref{thm:LeviPoisson},
there exists a constant $C$, such that for any positive number $\epsilon < 1/4$,
there is a natural number $K = K(\epsilon)$ and a positive number $\rho =
\rho(\epsilon)$, such that for any $l \geq K$ we can construct a local analytic
coordinate system $(x^l_1,...,y^l_{n-m})$ as in the previous sections, with the
following additional properties (using the previous notations) :

$(i)_l$ (Chains of balls) The ball $B_{l,exp(1/l)\rho}$ is well-defined, and if $l >
K$ we have
\begin{equation}
\label{eqn:inclusion0}  B_{l-1,exp(\frac{1}{l} - \frac{2\epsilon}{l^2})\rho} \subset
B_{l,exp(1/l)\rho} \subset B_{l-1,exp(\frac{1}{l} + \frac{2\epsilon}{l^2})\rho} \ .
\end{equation}

$(ii)_l$ (Norms of changes) If $l > K$ then we have
\begin{equation}
|\psi_l|_{l-1,exp(\frac{1}{l-1} - \frac{\epsilon}{(l-1)^2})\rho} < \rho \ .
\end{equation}

$(iii)_l$ (Norms of the Poisson structure) :
\begin{equation}
|\Pi|_{l,exp(1/l)\rho} \leq C . exp(-1/\sqrt{l})\rho \ .
\end{equation}

\end{prop}

Theorem \ref{thm:LeviPoisson} follows immediately from the first part of Proposition
\ref{prop:mainPoisson} and the following lemma:

\begin{lemma}
\label{lemma:convergence} If Condition $(i)_l$ of Proposition \ref{prop:mainPoisson}
is satisfied for all $l \geq K$ (where $K$ is some finite number), then the formal
coordinate system $(x^\infty_1,...,x^\infty_m,y^\infty_1,...,y^\infty_{n-m})$ is
convergent (i.e. locally analytic).
\end{lemma}

The main idea behind Lemma \ref{lemma:convergence} is that, if Condition $(i)_l$ is
true for any $l \geq K$, then the infinite intersection $\bigcap_{l=K}^\infty
B_{l,\exp(1/l)\rho}$ contains an open neighborhood of $0$, implying a positive
radius of convergence.

The second and third parts of Proposition \ref{prop:mainPoisson} are needed for the
proof of the first part. Proposition \ref{prop:mainPoisson} will be proved by
recurrence : By taking $\rho$ small enough, we can obviously achieve Conditions
$(iii)_K$ and $(i)_K$ (Condition $(ii)_K$ is void). Then provided that $K$ is large
enough, when $l \geq K$ we have that Condition $(ii)_l$ implies Conditions $(i)_l$
and $(iii)_l$, and Condition $(iii)_l$ in turn implies Condition $(ii)_{l+1}$.
%And if Condition $(i)_l$ is satisfied for all $l \geq K$, then the radius of convergence
%of the coordinate system $(x^\infty_1,...,y^\infty_{n-m})$ is positive.
In other words, Proposition \ref{prop:mainPoisson} follows directly from the
following three lemmas :

\begin{lemma}
\label{lemma:psi} There exists a finite number $K$ (depending on $\epsilon$) such
that if Condition $(iii)_{l}$ (of Proposition \ref{prop:mainPoisson}) is satisfied
and $l \geq K$ then Condition $(ii)_{l+1}$ is also satisfied.
\end{lemma}

\begin{lemma}
\label{lemma:balls} There exists a finite number $K$ (depending on $\epsilon$) such
that if Condition $(ii)_{l+1}$ is satisfied and $l \geq K$ then Condition
$(i)_{l+1}$ is also satisfied.
\end{lemma}

\begin{lemma}
\label{lemma:Pi} There exists a finite number $K$ (depending on $\epsilon$) such
that if Conditions $(ii)_{l+1}$ and $(iii)_l$ are satisfied and $l \geq K$ then
Condition $(iii)_{l+1}$ is also satisfied.
\end{lemma}

The lemmas of this section will be proved in detail in the subsequent section. Let
us mention here only the main ingredients behind the last three ones: The proof of
Lemma \ref{lemma:balls} and Lemma \ref{lemma:Pi} is straightforward and uses only
the first part of Lemma \ref{lemma:2norms}. Lemma \ref{lemma:psi} (the most
technical one) follows from the estimates on the primitives of cocycles as provided
by Lemma \ref{lemma:inverse2}.

\section{Proof of technical lemmas}

In this sections we will prove the lemmas stated in the previous section. \\

{\it Proof of Lemma \ref{lemma:2norms}}. Let $f$ be a local analytic function in
$(\bbC^n,0)$. To make an estimate on $df$, we use the Cauchy integral formula. For
$z \in B_\rho$, denote by $\gamma_i$ the following circle : $\gamma_i = \{ v \in
\bbC^n \ | \ v_j = z_j \ {\rm if} \ j \neq i \ , \ |v_i - z_i| = \epsilon\rho/2l^2
\}$. Then $\gamma_i \subset B_{(1+\epsilon/l^2)\rho}$, and we have
$$
\left| \frac{\partial f}{\partial z_i} (z) \right| = \frac{1}{2\pi} \left|
\oint_{\gamma_i} \frac{f(v) dv}{(v-z)^2}\right| \leq \frac{2l^2}{\epsilon\rho}
|f|_{(1+\epsilon/2l^2)\rho} \ ,
$$
which implies that $exp(2^{l/2}) |f|_{(1+\epsilon/2l^2)\rho} \geq {\rho}|df|$ when
$l$ is large enough.

Now let $f \in \cO_{l}$ such that $|f|_{(1+ \epsilon/l^2)\rho} < \infty$. We want to
show that if $x \in B_{(1 + \epsilon/2l^2)\rho}$ then $|f(x)| \leq exp(2^{l/2})
|f|_{(1+ \epsilon/l^2)\rho}$ (provided that $l$ is large enough compared to
$1/\epsilon$). Fix a point $x \in B_{(1 + \epsilon/2l^2)\rho}$ and consider the
following holomorphic function of one variable : $g(z) = f (\frac{x}{|x|}z)$. This
function is holomorphic in the complex 1-dimensional disk $B^1_{(1+
\epsilon/l^2)\rho}$ of radius $(1+ \epsilon/l^2)\rho$, and is bounded by $|f|_{(1+
\epsilon/l^2)\rho}$ in this disk. Because $f \in \cO_l$, we have that $g(z)$ is
divisible by $z^{2^l}$, that is $g(z) / z^{2^l}$ is holomorphic in $B^1_{(1+
\epsilon/l^2)\rho}$. By the maximum principle we have
$$\frac{|f(x)|}{|x|^{2^l}} = \left| \frac{g(|x|)}{|x|^{2^l}}\right| \leq
max_{|z| = (1 + \epsilon/l^2)\rho} \left| \frac{g(z)}{z^{2^l}} \right| \leq
\frac{|f|_{(1+ \epsilon/l^2)\rho}}{((1+\epsilon/l^2)\rho)^{2^l}} \ ,$$ which implies
that
$$
|f(x)| \leq (\frac{1 + \epsilon/2l^2}{1 + \epsilon/l^2})^{2^l}|f|_{(1+
\epsilon/l^2)\rho} \thickapprox exp(-\frac{2^l}{2\epsilon l^2}) |f|_{(1+
\epsilon/l^2)\rho} \leq exp(-2^{l/2}) |f|_{(1+ \epsilon/l^2)\rho}
$$
(when $l$ is large enough). Thus we have proved that there is a finite number $K$
depending on $\epsilon$ such that
$$
|f|_{(1+\epsilon/l^2)\rho} \geq exp(2^{l/2})|f|_{(1+ \epsilon/2l^2) \rho}
$$
for any $l > K$ and any $f \in \cO_l$.

To compare the norms of $f$, we use Cauchy-Schwartz inequality : for $f =
\sum_{\alpha \in \bbN^k} c_\alpha \prod_i z_i^{\alpha_i}$ and $|z| =
(1-\epsilon/2l^2)\rho$ we have
$$
\begin{array}{l}
|f(z)| \leq \sum_{\alpha \in \bbN^k} |c_\alpha| \prod_i |z_i|^{\alpha_i} \leq \cr
\leq \left( \sum_\alpha |c_\alpha|^2 \frac{\alpha ! (n-1) !}{(|\alpha| + n -1) !}
\rho^{2|\alpha|} \right)^{1/2} . \left( \sum_{\alpha} \frac{(|\alpha| + n -1)
!}{\alpha ! (n-1) !}  \rho^{-2|\alpha|} \prod_i |z_i|^{2\alpha} \right)^{1/2} = \cr
= \|f\|_\rho . \left( 1 - \sum_i \frac{|z_i|^2}{\rho^2} \right)^{-n/2} = \|f\|_\rho
. (1 - (1 - \epsilon/2l^2)^2)^{-n/2} \leq \frac{(2l)^n}{\epsilon^{n/2}} \|f\|_{\rho}
\ .
\end{array}
$$
It means that for any local analytic function $f$ we have
\begin{equation}
|f|_{(1-\epsilon/2l^2)\rho} \leq \frac{(2l)^n}{\epsilon^{n/2}} \|f\|_{\rho} \ .
\end{equation}
Now if $f \in \cO_l$, we can apply Inequality (\ref{eqn:changeradius}) to get
$$
|f|_{(1-\epsilon/l^2)\rho} \leq exp(-2^{l/2}) |f|_{(1-\epsilon/2l^2)\rho} \leq
\frac{(2l)^n}{\epsilon^{n/2}} exp(-2^{l/2}) \|f\|_{\rho} \leq \|f\|_{\rho} \ ,
$$
provided that $l$ is large enough compared to $1/\epsilon$.
Lemma \ref{lemma:2norms} is proved. \hfill $\square$ \\

{\it Proof of Lemma \ref{lemma:convergence}}. The main point is to show that the
limit $\bigcap_{l=K}^{\infty} B_{l,\rho}$ contains a ball $B_r$ of positive radius
centered at $0$. Then for $x \in B_r$, we have $x \in B_{l,\rho}$, implying
$\|(z^l_1(x),...,z^l_{n}(x))\| < \rho$ is uniformly bounded, which in terms implies
that the formal functions $z^\infty_i = \lim_{l \to \infty} z^l_i$ are analytic
functions inside $B_r$ (recall that $(z^l_1,...,z^l_n)$ is obtained from
$(x^l_1,...,y^l_{n-m})$ by a constant linear transformation $(A_{ij})$ which does
not depend on $l$).

Recall the following fact of complex analysis, which is a consequence of the maximum
principle : if $g$ is a complex analytic map from a complex ball of radius $\rho$ to
some linear Hermitian space such that $g(0) = 0$ and $|g(x)| \leq C$ for all $|x| <
\rho$ and some constant $C$, then we have $|g(x)| \leq C|x|/\rho$ for all $x$ such
that $|x| < \rho$. If $l_1,l_2 \in \mathbb N$ and $r_1,r_2 > 0, s > 1$, then
applying this fact we get :

\begin{equation}
\label{eqn:maximal} {\rm If} \ B_{l_1,r_1} \subset B_{l_2,r_2} \ {\rm then} \
B_{l_1,r_1/s} \subset B_{l_2,r_2/s} \ .
\end{equation}

(Here $r_1$ plays the role of $\rho$, $r_2$ plays the role of $C$, and the
coordinate transformation from $(z^{l_1}_1,...,z^{l_1}_{n})$ to
$(z^{l_2}_1,...,z^{l_2}_{n})$ plays the role of $g$ in the previous statement).

Using Formula (\ref{eqn:maximal}) and Condition $(i)_l$ recursively, we get

\begin{equation}
B_{l,\rho} \supset B_{l-1,exp(-1/l^2) \rho} \supset B_{l-2,exp(-1/l^2 -
1/(l-1)^2)\rho} \supset ... \supset B_{K,exp(-\sum_{k=K}^l 1/k^2)\rho} \ .
\end{equation}

Since $c = exp(-\sum_{k=K}^\infty 1/k^2)$ is a positive number, we have
$\bigcap_{l=K}^{\infty} B_{l,\rho} \supset B_{K,c\rho}$, which clearly contains an
open
neighborhood of $0$. Lemma \ref{lemma:convergence} is proved. \hfill $\square$ \\

{\it Proof of Lemma \ref{lemma:balls}}. Suppose that Condition $(ii)_{l+1}$ is
satisfied. For simplicity of exposition, we will assume that the coordinate system
$(z^l_1,...,z^l_n)$ coincides with the coordinate system $(x^l_1,...,y^l_{n-m})$
(The more general case, when $(z^l_1,...,z^l_n)$ is obtained from
$(x^l_1,...,y^l_{n-m})$ by a constant linear transformation, is essentially the
same). Suppose that we have
\begin{equation}
|\psi_{l+1}|_{l,exp (1/l - \epsilon/l^2)\rho} < \rho \ .
\end{equation}
Then it follows from  Lemma \ref{lemma:2norms} that, provided that $l$ is large
enough :
\begin{equation}
|d\psi_{l+1}|_{l,exp (1/l - 2\epsilon/l^2)\rho} < 1/2n \ .
\end{equation}
(In order to define $|d\psi_{l+1}|_{l,exp (1/l - 2\epsilon/l^2)\rho}$, consider
$d\psi_{l+1}$ as an $n^2$-vector valued function in variables
$(z_1^{l},...,z_{n}^{l})$). Hence the map $\phi_{l+1} = Id + \psi_{l+1}$ is
injective in $B_{l,exp (1/l - 2\epsilon/l^2)\rho}$ : if $x, y \in B_{l,\rho_l}, x
\neq y$, then $\|\phi_{l+1}(x) - \phi_{l+1}(y)\| \geq \|x-y\| - \|\psi_{l+1}(x) -
\psi_{l+1}(y)\| \geq \|x-y\| - n|d\psi_{l+1}|_{exp (1/l - 2\epsilon/l^2)\rho}
\|x-y\| \geq (1 - 1/2)\|x-y\| > 0$.  (Here $(x-y)$ means the vector $(z^l_1 (x) -
z^l_1(y),..., z^l_n(x) - z^l_n(y))$, i.e. their difference is taken with respect to
the coordinate system $(z_1^{l},...,z_{n}^{l})$).

It follows from Lemma \ref{lemma:2norms} that $|\phi_{l+1}|_{l,exp (1/l -
2\epsilon/l^2)\rho} = |Id + \psi_{l+1}|_{l,exp (1/l - 2\epsilon/l^2)\rho} \leq
|Id|_{l,exp (1/l - 2\epsilon/l^2)\rho} + |\psi_{l+1}|_{l,exp (1/l -
2\epsilon/l^2)\rho} < exp (1/l - 2\epsilon/l^2)\rho + \frac{\epsilon}{4l^2} exp (1/l
- 2\epsilon/l^2)\rho < exp(1/l - \epsilon/l^2) \rho$. In other words, we have
\begin{equation}
\label{eqn:inclusion1} \phi_{l+1} (B_{l,exp (1/l - 2\epsilon/l^2)\rho}) \subset
B_{l,exp(1/l - \epsilon/l^2) \rho} \ .
\end{equation}

Applying Formula (\ref{eqn:maximal}) to the above relation, noticing that $1/l -
2\epsilon/l^2 > 1/(l+1)$, and simplifying the obtained formula a little bit, we get
\begin{equation}
\label{eqn:inclusion2} \phi_{l+1} (B_{l,exp(1/(l+1)- 2\epsilon/(l+1)^2)\rho})
\subset B_{l,exp(1/(l+1))\rho} \ .
\end{equation}

We will show that $\phi_{l+1}^{-1}$ is well-defined in $B_{l,exp(1/(l+1))\rho}$, and
\begin{equation}
\label{eqn:inclusion3} \phi_{l+1}^{-1} (B_{l,exp(1/(l+1))\rho}) =
B_{l+1,exp(1/(l+1))\rho} \subset B_{l,exp(1/(l+1)+ 2\epsilon/(l+1)^2 ) \rho} \ .
\end{equation}
Indeed, if we denote by $S_{l,exp (1/l - 2\epsilon/l^2)\rho}$ the boundary of
$B_{l,exp (1/l - 2\epsilon/l^2)\rho}$, then $\phi_{l+1}(S_{l,exp (1/l -
2\epsilon/l^2)\rho})$ lies in $B_{l,exp (1/l - \epsilon/l^2)\rho}$ and is homotopic
to $S_{l,exp (1/l - 2\epsilon/l^2)\rho}$ via a homotopy which does not intersect
$B_{l,exp(1/(l+1))\rho}$. It implies (via the classical Brower's fixed point
theorem) that $\phi_{l+1} (B_{l,exp (1/l - 2\epsilon/l^2)\rho})$ must contain
$B_{l,exp(1/(l+1))\rho}$. Because $\phi_{l+1}$ is injective in $(B_{l,exp (1/l -
2\epsilon/l^2)\rho})$, it means that the inverse map is well-defined in
$B_{l,exp(1/(l+1))\rho}$, with $ \phi_{l+1}^{-1}(B_{l,exp(1/(l+1))\rho}) \subset
B_{l,exp (1/l - 2\epsilon/l^2)\rho}$. In particular, $B_{l+1,exp(1/(l+1))\rho} =
\phi_{l+1}^{-1}(B_{l,exp(1/(l+1))\rho})$ is well-defined. Lemma \ref{lemma:balls}
then follows from (\ref{eqn:inclusion2}) and (\ref{eqn:inclusion3}).
\hfill $\square$ \\

{\it Proof of Lemma \ref{lemma:psi}}. Suppose that Condition $(iii)_l$ is satisfied.
Then according to (\ref{eqn:fij}) we have :
\begin{equation}
\label{eqn:estf}
\begin{array}{c}
\|f^l_{ij}\|_{l,exp(1/l)\rho} \leq |f^l_{ij}|_{l,\exp(1/l)\rho} = |\{x^l_i,x^l_j\} -
\sum_k c_{ij}^k x^l_k|_{l,exp(1/l)\rho} \leq \cr \leq C_1|\Pi|_{l,exp(1/l)\rho} +
\sum_k |c_{ij}^k\|x^l_k|_{l,\rho} \leq C_1 . C . \rho + C_2 . exp(1/l) \rho \sum_k
|c_{ij}^k| < C_3 \rho \ ,
\end{array}
\end{equation}
where $C_3$ is some positive constant (which does not depend on $l$).

We can apply the above inequality $\|f^l_{ij}\|_{l,exp(1/l)\rho} < C_3 \rho$ and
Lemma \ref{lemma:inverse2} to find a positive constant $C_4$ (which does not depend
on $l$) and a solution $w^l_i$ of (\ref{eqn:wi2}), such that
\begin{equation}
\label{eqn:estw} \|w^l_i\|_{l,exp(1/l)\rho} < C_4 \rho \ .
\end{equation}

Together with Lemma \ref{lemma:2norms}, the above inequality yields

\begin{equation}
\label{eqn:estdw} |dw^l_i|_{l,exp(1/l - \epsilon/2l^2)\rho} < C_4 ,
\end{equation}
provided that $l$ is large enough. Applying Lemma \ref{lemma:2norms} and the
assumption that $|\Pi|_{l,exp(1/l)\rho} < C \rho$ to the above inequality, we get

\begin{equation}
|\{w_i^l,y_j^l\}|_{l,exp(1/l - \epsilon/2l^2)\rho} < C_5 \rho
\end{equation}
for some constant $C_5$ (which does not depend on $l$). Using this inequality, and
inequalities similar to (\ref{eqn:estf}), we get that the norm $\|.\|_{l,exp(1/l -
\epsilon/2l^2)\rho}$ of the 1-cocycle given in Formula (\ref{eqn:1cocycle}) is
bounded from above by $C_6 \rho$, where $C_6$ is some constant which does not depend
on $L$. Using Lemma \ref{lemma:inverse2}, we find a solution $v_i^L$ to Equation
\ref{eqn:vi} such that
\begin{equation}
\label{eqn:estv} \|v^l_i\|_{l,exp(1/l - \epsilon/2l^2)\rho} < C_6 \rho \ ,
\end{equation}
where $C_6$ is some constant which does not depend on $l$. Lemma \ref{lemma:psi} (fr
$l$ large enough compared to $C_6$) now follows directly from Inequalities
(\ref{eqn:estw}),
(\ref{eqn:estv}) and Lemma \ref{lemma:2norms}. \hfill $\square$ \\

{\it Proof of Lemma \ref{lemma:Pi}}. Suppose that Condition $(ii)_{l+1}$ is
satisfied. By Lemma \ref{lemma:balls}, Condition $(i)_{l+1}$ is also satisfied. In
particular,
$$B_{l+1,exp(1/(l+1))\rho} \subset B_{l,exp(1/(l+1) + 2\epsilon/(l+1)^2)\rho}
\subset  B_{l,exp(1/l - 2\epsilon/l^2)\rho}  $$ (for $\epsilon < 1/4$ and $l$ large
enough). Thus we have
\begin{equation}
|\{z^{l+1}_i,z^{l+1}_j\}|_{l+1,exp(1/(l+1))\rho} \leq
|\{z^{l+1}_i,z^{l+1}_j\}|_{l,exp(1/l - 2\epsilon/l^2)\rho} \leq T^1 + T^2 + T^3 +
T^4
\end{equation}
where
\begin{equation}
\begin{array}{l}
T^1 = |\{z^{l}_i,z^{l}_j\}|_{l,exp(1/l - 2\epsilon/l^2)\rho} \ , \\
T^2 = |\{z^{l+1}_i - z^l_i,z^{l+1}_j\}|_{l,exp(1/l - 2\epsilon/l^2)\rho} \ , \\
T^3 = |\{z^{l+1}_i,z^{l+1}_j-z^l_j\}|_{l,exp(1/l - 2\epsilon/l^2)\rho} \ , \\
T^4 = |\{z^{l+1}_i - z^l_i,z^{l+1}_j- z^l_j\}|_{l,exp(1/l - 2\epsilon/l^2)\rho} \ .
\end{array}
\end{equation}

For the first term, we have
$$
T^1 \leq |\{z^{l}_i,z^{l}_j\}|_{l,exp(1/l)\rho} \leq |\Pi|_{l,exp(1/l)\rho} \leq C.
exp(-1/\sqrt{l})\rho \ .
$$
Notice that $C exp(-1/\sqrt{l+1}) \rho - C exp(-1/\sqrt{l})\rho > \frac{C}{l^2}
\rho$ (for $l$ large enough). So to verify Condition $(iii)_{l+1}$, it suffices to
show that $T^2 + T^3 + T^4 < \frac{C}{l^2} \rho$. But this last inequality can be
achieved easily (provided that $l$ is large enough) by Conditions $(ii)_{l+1}$,
$(iii)_{l}$ and Lemma \ref{lemma:2norms}. Lemma \ref{lemma:Pi} is proved.
\hfill $\square$ \\

\section{Lie algebroids}

Let $A = ( \bbK^N \times (\bbK^n,0) \rightarrow (\bbK^n,0)\ , \ [\ ,\ ] \ , \ \#)$
be a local analytic Lie algebroid, with Lie bracket $[\ ,\ ]$ and anchor map $\#$.
It is well-known that (see e.g. \cite{CaWe}), on the total space of the dual bundle
$A^\ast = (\bbK^N)^\ast \times (\bbK^n,0) \rightarrow (\bbK^n,0)$ there is a unique
natural Poisson structure associated to $A$ (called the dual Lie-Poisson structure),
defined as follows. By duality, consider sections of $A$ as fiber-wise linear
functions on (the total space of) $A^\ast$. Let $(x_1,...,x_m)$ be a coordinate
system of $(\bbK^n,0)$, and $(s_1,...,s_N)$ be a basis of the space of sections of
$A$. Then $(x_1,...,x_m,s_1,...,s_N)$ is a coordinate system for $A^\ast$, and the
Poisson bracket on $A^\ast$ is given by the following formula :
\begin{equation}
\begin{array}{l}
\{s_i,s_j\} = [s_i,s_j] \ ,  \\
\{s_i,x_j\} = \#s_i (x_j) \ , \\
\{x_i,x_j\} = 0 \ .
\end{array}
\end{equation}
The above Poisson structure is fiber-wise linear in the sense that the Poisson
bracket of two fiber-wise linear functions is again a fiber-wise linear function,
the Poisson bracket of a fiber-wise linear function and a base function is a base
function, and the Poisson bracket of two base functions is zero. Conversely, it is
clear that any such a Poisson structure on a bundle $A^\ast = (\bbK^N)^\ast \times
(\bbK^n,0) \rightarrow (\bbK^n,0)$ corresponds to a Lie algebroid structure on the
dual bundle $\bbK^N \times (\bbK^n,0) \rightarrow (\bbK^n,0)$.

It is easy to see that, to prove Theorem \ref{thm:LeviAlgebroid}, it suffices to
find a Levi factor for the dual Lie-Poisson structure, which consists of fiber-wise
linear functions. The existence of a Levi factor for the Poisson structure on
$A^\ast$ is provided by Theorem \ref{thm:LeviPoisson}. We only have to make sure
that this Levi factor can be chosen so that it consists of fiber-wise linear
functions. In order to see it, one makes the following modifications to the
construction of Levi decomposition given in Section \ref{section:formal} :

- After Step $l$ ($l \geq 0$), we will get a local coordinate system
$$(s^l_1,...,s^l_m,v^l_1,...,v^l_{N-m},x^l_1,...,x^l_n)$$
of $A^\ast$ with the following properties : $x^l_1,...,x^l_n$ are base functions
(i.e. functions on $(\bbK^n,0)$ ; $s^l_1,...,s^l_m,v^l_1,...,v^l_{N-m}$ are
fiber-wise linear functions (i.e. they are sections of $A$) ; $\{s^l_i,s^l_j\} -
\sum_k c_{ij}^k s^l_k = O(|x|^{2^l})$ ; $\{s^l_i,v^l_j\} - \sum_k a_{ij}^k v^l_k =
O(|x|^{2^l})$ ; $\{s^l_i,x^l_j\} - \sum_k b_{ij}^k x^l_k = O(|x|^{2^l+1})$ . Here
$c_{ij}^k, a_{ij}^k, b_{ij}^k$ are structural constants as appeared in the statement
of Theorem \ref{thm:LeviAlgebroid}.

 - Replace the space $\cO$ of all local analytic functions by the subspace of local analytic
functions which are fiber-wise linear. Similarly, replace the space $\cO_l$ of local
analytic functions without terms of order $\leq 2^l$ by the subspace of fiber-wise
linear analytic functions without terms of order $\leq 2^l$.

- Replace $\cY^l$ by the subspace of vector fields of the following form :
$$
\sum_{i=1}^{N-m} p_i \partial / \partial v^l_i +  \sum_{i=1}^{n} q_i \partial /
\partial x^l_i \ ,
$$
where $p_i$ are fiber-wise linear functions and $q_i$ are base functions. For the
replacement of $\cY^l_k$, we require that $p_i$ do not contains terms of order $\leq
2^k - 1$ in variables $(x_1,...,x_n)$, and $q_i$ do not contains terms of order
$\leq 2^k$.

One checks that the above subspaces are invariant under the $\fg$-actions introduced
in Section \ref{section:formal}, and the cocycles introduced there will also live in
the corresponding quotient spaces of these subspaces. Details are left to the
reader. \hfill $\square$ \\

The smooth version of the main results of this paper is considered in a separate
work in collaboration with Philippe Monnier \cite{MoZu-Levi2002}. The results of
\cite{MoZu-Levi2002} generalize Conn's smooth linearization theorem for smooth
Poisson structures with a compact semisimple linear part \cite{Conn-Smooth1985}, and
imply the local smooth linearizability of smooth Lie algebroids with a compact
semisimple linear part.
\\

{\bf Acknowledgements}. I would like to thank  Jean-Paul Dufour for many discussions
on the subject of this paper, and for providing me with various historical remarks.
I'm also thankful to the referee for numerous suggestions which helped improve the
presentation of this paper.

\bibliographystyle{amsplain}

\begin{thebibliography}{ABC}

\bibitem{Arnold-Small1963}
V.I. Arnold, {\it Small denominators and problems of stability motion in classical
and celestial mechanics}, Russian Math. Surveys, 18 (1963), 85-191.

\bibitem{CaWe}
A. Cannas da Silva, A. Weinstein, Geometric models for noncommutative algebras,
Berkeley Mathematics Lectures, Vol. 10, AMS, 1999.

\bibitem{Cerveau}
D. Cerveau, {\it Distributions involutives singulières.} Ann. Inst. Fourier
(Grenoble) 29 (1979), no. 3, xii, 261--294.

\bibitem{Conn}
J.F. Conn, {\it Normal forms for analytic Poisson structures.} Ann. of Math. (2) 119
(1984), no. 3, 577--601; and {\it Correction to: ``Normal forms for analytic Poisson
structures''}, Ann. of Math. (2) 125 (1987), no. 2, 433--436.

\bibitem{Conn-Smooth1985}
J.F. Conn, {\it Normal forms for smooth Poisson structures.} Ann. of Math. (2) 121
(1984), 565-593.

\bibitem{Dufour-algebroid}
J.-P. Dufour, {\it Normal forms of Lie algebroids}, Banach Center Publications, Vo.
54 (2001), 35-41.

\bibitem{DuZu-affn2002}
J.-P. Dufour, Nguyen Tien Zung, {\it Nondegeneracy of the Lie algebra
$\mathfrak{aff}(n)$}, to appear in Comptes Rendus Acad. Sci. Paris.

\bibitem{Fernandes-Algebroid2001}
L. R. Fernandes, {\it Lie algebroids, holonomy and characteristic classes}, Adv. in
Math., 170 (2002), 119-179.

\bibitem{GuSt-Linearization1968}
V. Guillemin, S. Sternberg, {\it Remarks on a paper of Hermann}, Trans. Amer. Math.
Soc., 130 (1968), 110-116.

\bibitem{Jacobson}
N. Jacobson, Lie algebras, Interscience Publishers, New York, 1962.

\bibitem{Molinier}
J.-C. Molinier, {\it Linéarisation de structures de Poisson}, Ph.D. thesis,
Montpellier 1993.

\bibitem{MoZu-Levi2002}
Ph. Monnier, Nguyen Tien Zung, {\it Levi decomposition for smooth Poisson
structures}, preprint math.DG/0209004.

\bibitem{Wade}
A. Wade, {\it  Normalisation formelle de structures de Poisson.} C. R. Acad. Sci.
Paris S\'er. I Math. {\bf 324} (1997), 531--536.

\bibitem{Weinstein-local}
A. Weinstein, {\it The local structure of Poisson manifolds.} J. Differential Geom.
18 (1983), no. 3, 523--557.

\bibitem{Weinstein-algebroid}
A. Weinstein, {\it  Linearization problems for Lie algebroids and Lie groupoids.}
Conference Moshé Flato 1999 (Dijon). Lett. Math. Phys. 52 (2000), no. 1, 93--102.

\bibitem{Zung-Birkhoff2001}
Nguyen Tien Zung, {\it Convergence versus integrability in Birkhoff normal form},
preprint math.DS/0104279 (2001).

\bibitem{Zung-PD2002}
Nguyen Tien Zung {\it Convergence versus integrability in Poincaré-Dulac normal
form}, Math. Res. Lett. 9 (2002), 217-228.

\end{thebibliography}

\end{document}